\documentclass[12pt]{amsart}
\newif\ifpdf
    \ifx\pdfoutput\undefined
    \pdffalse 
    \else
    \pdfoutput=1 
    \pdftrue
    \fi

    \ifpdf
    \usepackage[pdftex]{graphicx}
    \else
    \usepackage{graphicx}
    \fi
\usepackage{amsmath,amssymb,amsthm,epsfig}

\newtheorem{theorem}{Theorem}[section]
\newtheorem{proposition}[theorem]{Proposition}
\newtheorem{lemma}[theorem]{Lemma}
\newtheorem{corollary}[theorem]{Corollary}
\newtheorem{example}[theorem]{Example}

\input{epsf.tex}

\newcommand\F{\mathcal F}

\newcommand\x{x_0}
\newcommand\inv{x_0^{-1}}
\newcommand\Ll{L_L}
\newcommand\rni{R_{NI}}

\begin{document}
\ifpdf
    \DeclareGraphicsExtensions{.pdf, .jpg, .tif}
    \else
    \DeclareGraphicsExtensions{.eps, .jpg}
    \fi

\title{Combinatorial properties of Thompson's group $F$}
\author{Sean Cleary and Jennifer Taback}
\thanks{The first author acknowledges support from PSC-CUNY grant \#63438-0032.
The second author would like to thank the University of
Utah for their hospitality during the writing of this paper.}

\begin{abstract}
We study some combinatorial properties of the word metric
of  Thompson's group $F$ in the standard two generator finite presentation.
We explore connections between the tree pair diagram
representing an element $w$ of $F$, its normal form in the infinite presentation,
its word length, and minimal length representatives of it.
We estimate word length in terms of the number and type
of carets in the tree pair diagram and show sharpness of those
estimates.  In addition we explore some properties of the Cayley
graph of $F$ with respect to the two generator  finite presentation.
Namely, we exhibit the form of ``dead end'' elements in this Cayley
graph, and show that it has no ``deep pockets''.
Finally, we discuss a simple method for constructing minimal
length representatives for strictly positive or negative
words.
\end{abstract}

\maketitle

\section{Introduction}
\label{sec:intro}

Thompson's group $F$ has been studied extensively in many
different branches of mathematics, including group theory,
dynamics, homotopy theory and logic.  Algebraically, it is most commonly
understood in two different forms: via a finite presentation and
an infinite presentation.  The infinite presentation ${\mathcal
P}$ has  simple relators which are
conveniently  manipulated and understood, as well as a unique normal form for elements.
The finite presentation ${\mathcal F}$ has two generators and two relators, but the
relators  are more complicated and there is no longer a convenient
set of normal forms.  The elements of $F$ can also be interpreted as pairs of
finite binary rooted trees with the same number of carets.

In this paper we discuss many interesting combinatorial properties of
Thompson's group $F$.  These properties are derived from  the
relationship between the normal form of elements of $F$ and the pairs
 of finite binary rooted trees used to represent elements of $F$.
The combinatorial properties we describe have applications to
estimating word length in $F$, counting and determining caret
types and lead to algorithms for constructing minimal length
paths. The sections of this paper are organized as follows.

\begin{itemize}
\item
In {\bf Section 2} we present brief introduction to Thompson's group $F$,
including Fordham's method of calculating word length \cite{blake}.  We
detail the bijective process which transforms a tree pair diagram
representing an element of $F$ into its unique normal form.

\item
In {\bf Section 3} we apply Fordham's method to immediately obtain an
estimate of the word length $|w|$ of an element $w \in F$ from the
tree pair diagram representing $w$ in the word metric arising from the
finite presentation $\F$.  We give examples which show that
the constants in the estimate are sharp.

\item
In {\bf Section 4} we use Fordham's method of calculating word length to
explore an interesting phenomenon which occurs in the Cayley graph of
$F$ with respect to the standard  two generator finite presentation.
Namely, there are {\em dead end elements} $w$ with the property that
$|w \alpha| = |w| - 1$ for all generators $\alpha \in \{ \x^{\pm 1}, \ x_1^{\pm 1} \}$, where
$|w|$ denotes word length.  Fordham
\cite{blake} remarks that some of these dead end elements have a
particular form; we give a general form for all dead end elements in $F$
and describe some limits to stronger forms of this behavior, called ``deep pockets.''

\item
In {\bf Section 5} describe the combinatorial relationship between the
normal form of an element $w \in F$ and the number and types of carets
in the tree pair diagram representing $w$.

\item
In {\bf Section 6} we present a method of
constructing minimal length paths in the standard two generator presentation
for strictly positive or negative words in
$F$.
\end{itemize}

\section{Thompson's group $F$}
\label{sec:F}

Thompson's group is best understood combinatorially using the two
presentations mentioned above, the finite presentation
$$
{\mathcal F} = \langle x_0,x_1 |
[x_0x_1^{-1},x_0^{-1}x_1x_0],[x_0x_1^{-1},x_0^{-2}x_1x_0^2]
\rangle$$ and the infinite presentation
$$
{\mathcal P} = \langle x_k, \ k \geq 0 | x_i^{-1}x_jx_i = x_{j+1}
\ \text{ if }i<j \rangle.$$

A convenient set of normal forms for elements of $F$ in
the infinite presentation ${\mathcal P}$ is given by $x_{i_1}^{r_1}
x_{i_2}^{r_2}\ldots x_{i_k}^{r_k} x_{j_l}^{-s_l} \ldots
x_{j_2}^{-s_2} x_{j_1}^{-s_1} $ where $r_i, s_i >0$, $i_1<i_2
\ldots < i_k$ and $j_1<j_2 \ldots < j_l$.   To obtain a unique normal
form for each element, we add the condition that when
 both $x_i$ and $x_i^{-1}$
occur, so does $x_{i+1}$ or  $x_{i+1}^{-1}$, as discussed by
Brown and Geoghegan
\cite{bg:thomp}.   We will always mean unique normal form
when we refer to a word $w$ in normal form.

Analytically, we can regard $F$ as the group of
orientation-preserving piecewise-linear homeomorphisms from
$[0,1]$ to itself where each homeomorphism has only finitely many
singularities of slope, all such singularities lie in the dyadic
rationals ${\bf Z}[\frac12]$, and, away from the singularities,
the slopes are powers of $2$.

\subsection{Tree pair diagrams}
An element of $F$ can be interpreted geometrically via a tree pair
diagram,  which is a pair of rooted binary trees
$(T_-,T_+)$, each with the same number of exposed leaves, as described in
Cannon, Floyd and Parry \cite{cfp}.
An {\em exposed} leaf ends in a vertex of valence $1$, and we number
these exposed leaves from left to right, beginning with $0$.
We refer to a node together with the two
downward-directed edges from the node as a {\em caret}.
A caret $C$ may have a {\em right child},  a caret $C_R$ which is attached to the
right edge of $C$.  We can similarly define the {\em left child} $C_L$ of the
caret $C$.   The set of all carets which stem from the right leaf of a caret
$C$ is called the right subtree of $C$, and we can analogously define the
left subtree of $C$.

\begin{figure}\includegraphics[width=3in]{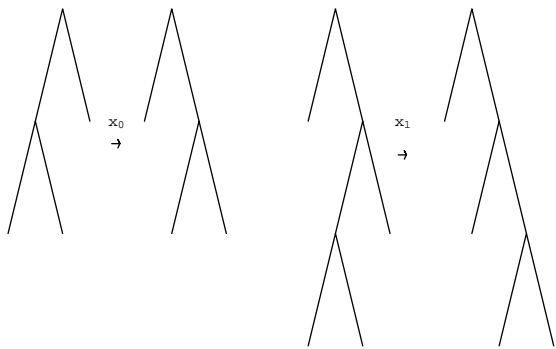}\\
\caption{The tree pair diagrams for the generators $\x$ and $x_1$ of
$\F$.
\label{fig:x0x1trees}}
\end{figure}

In a tree pair $(T_-,T_+)$, the tree
$T_-$ is called the {\em negative tree} and $T_+$ the {\em positive tree}.
This terminology is explained further in \S \ref{sec:exponents} below.
The equivalence between tree pair diagrams and
homeomorphisms of $[0,1]$ is described in \cite{cfp}.
In figure \ref{fig:x0x1trees} we give the tree pair diagrams for
the generators $\x$ and $x_1$ of $\F$.  In \S \ref{sec:exponents}
below the correspondence between the trees and the elements is explained.

A tree pair diagram is {\em unreduced} if both $T_-$ and $T_+$
contain a caret with two exposed leaves numbered $m$ and $m+1$.
There are many tree pair diagrams representing the same element of
$F$ but each element has a unique reduced tree pair diagram representing it.
When we write $(T_-,T_+)$ to represent an element
of $F$, we are assuming that the tree pair is reduced.

We refer the reader to Cannon, Floyd and Parry \cite{cfp} for an excellent introduction to
Thompson's group $F$, and to Cleary and Taback \cite{ct} for more details on
understanding geometrically the elements of $F$ as reduced tree pairs.   All of
the geometric facts used below are justified in \cite{ct}.

\subsection{Exponents in tree pair diagrams}
\label{sec:exponents}
There is a bijective correspondence between the tree pair diagram of $w=(T_-,T_+)$
and the normal form of $w$, described in \cite{cfp}.  In the tree pair $(T_-,T_+)$,
number the exposed leaves of $T_-$ and $T_+$ from left to right, beginning with $0$.
The {\em exponent} of the leaf labelled $k$, written $E(k)$, is
defined as the length of the maximal path consisting entirely of
left edges from $k$ which does not reach the right side of the
tree. Note that $E(k)=0$ for an exposed leaf labelled $k$ which is a right leaf of a
caret, as there is no path consisting entirely of left edges
originating from $k$.  In figure \ref{fig:x0x1trees}, number the exposed
leaves of the trees in the pair representing $x_1$ from left to
right, beginning with $0$.  Then the exponents of the leaves of $T_-$ are
all $0$, and the exponents of the leaves of $T_+$ are $0,1,0,0$, in order.
We refer the reader to \cite{ct} for a more detailed example of computing
exponents in a tree.

Once the exponents of the leaves in $T_-$ and $T_+$ have been computed, the
normal form of the element $w = (T_-,T_+)$ is easily obtained.  The positive
part of the normal form of $w$ is
$$
x_0^{E(0)} x_1^{E(1)} \cdots x_m^{E(m)}
$$
where $m$ is the number of exposed leaves in either tree, and the exponents are
obtained from the leaves of $T_+$.  The negative part of the normal form of $w$ is
similarly found to be
$$
x_m^{-E(m)} x_{m-1}^{-E(m-1)} \cdots x_0^{-E(0)}
$$
where the exponents are now computed from the leaves of $T_-$.  Note
that many of the exponents in the normal form as given above may be zero.

Similarly, given an element $x$ in normal form with respect to the
infinite generating set, it is possible to construct a tree pair
diagram  ($T_-, T_+$) so that each leaf has the correct exponent.
If $R$ is a right caret with a single exposed left leaf labelled $k$, then
$E(k) = 0$ by definition.  Thus, arbitrarily many right carets with no left subtrees
can be added to either $T_-$ or $T_+$ without affecting the normal
form to ensure that both trees have the same number of carets, and thus equivalently
the same number of exposed leaves.

\subsection{Fordham's method of calculating word length}
\label{sec:Fordham}

For an element $w$ of $F$, we let $|w|$ denote the word length of $w$ with
respect to the word metric arising from the finite presentation
$\F$.   Fordham \cite{blake} presents a method of calculating $|w|$ based
on the trees $T_-$ and $T_+$ in the tree pair diagram representing $w$.
He defines seven types of
carets that can be found in a rooted binary tree, and an
intricate system of weights assigned to different pairs of
caret types, which sum to $|w|$.  A detailed example of
calculating $|w|$ in this way can be found in \cite{ct}.

Let $T$ be a finite rooted binary tree.  The {\em left side} of
$T$ is the maximal path of left edges beginning at the root of
$T$. Similarly, we have the {\em right side} of $T$.  A caret in $T$
is a {\em left caret} if its left edge is on the left side of the
tree, a {\em right caret} if it is not the root and its right edge
is on the right side of the tree, and an {\em interior caret}
otherwise.  The carets and the exposed leaves of $T$ are
numbered independently, according to different methods.
As above, the exposed leaves are numbered from left to
right, beginning with $0$.
The carets in $T$ are numbered according to the infix ordering of nodes.
Caret $0$ is a left caret with an exposed left leaf
numbered $0$ in the leaf numbering.
According to the infix scheme, we number the left children of a caret before the caret
itself, and number the right children after numbering the caret.

\begin{figure}\includegraphics[width=3in]{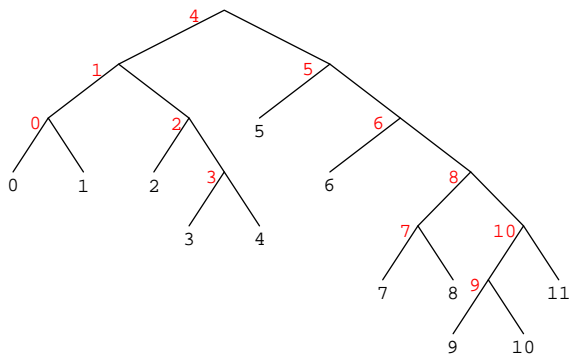}\\
\caption{A tree whose leaves are numbered left to right and whose nodes are numbered according to the
infix method.
\label{fig:infix}}
\end{figure}

In figure \ref{fig:infix} we give an example of a tree whose carets are
numbered according to the infix numbering method.
More examples of trees whose carets are numbered in this way can be
found in \cite{ct}.

Fordham classifies carets into seven disjoint types, as follows:
\begin{enumerate}

\item
$L_0$.  The first caret on the left side of the tree, with
caret number $0$. Every tree has exactly one caret of type $L_0$.

\item
$\Ll$.  Any left caret other than the one numbered $0$.

\item
$I_0$. An interior caret which has no right child.

\item
$I_R$. An interior caret which has a right child.

\item
$R_I$. Any right caret numbered $k$ with the property that caret
$k+1$ is an interior caret.

\item
$\rni$. A right caret which is not an $R_I$ but for which there is a
higher numbered interior caret.

\item
$R_0$. A right caret with no higher-numbered interior carets.

\end{enumerate}

The root caret is always considered to be a left caret and will
be of type $\Ll$ unless
it has no left children, in which case it would be the single caret
of type $L_0$.

The main result of  Fordham \cite{blake} is that the word length $|w|$
of $w = (T_-,T_+)$ can be computed from knowing the caret types of
the carets in the two trees, as long as they form a reduced pair,
 via the following process.  We number
the $k+1$ carets according to the infix method described above,
and for each $i$ with $0 \leq i \leq k$ we form the pair of caret
types consisting of the type of caret number $i$ in $T_-$ and the
type of caret number $i$ in $T_+$. The single caret of type $L_0$
in  $T_-$ will be paired with the single caret of type $L_0$ in
$T_+$, and for that pairing we assign a weight of 0. For all other
caret pairings, we assign weights according to the following
table.

\begin{center}

\begin{tabular}{|c|c|c|c|c|c|c|}

\hline
 & $R_0$ & $\rni$ & $R_I$ & $\Ll$ & $I_0$ & $I_R$ \\
 \hline

 $R_0$ & 0 & 2 & 2 & 1 & 1 & 3 \\ \hline
 $\rni$ & 2 & 2 & 2 & 1 & 1 & 3 \\ \hline
 $R_I$ & 2 & 2 & 2 & 1 & 3 & 3 \\ \hline
 $\Ll$ & 1 & 1 & 1 & 2 & 2 & 2 \\ \hline
 $I_0$ & 1 & 1 & 3 & 2 & 2 & 4 \\ \hline
$I_R$ & 3 & 3 & 3 & 2 & 4 & 4 \\ \hline
\end{tabular}

\end{center}

The main result of Fordham \cite{blake} is the following theorem.

\begin{theorem}[Fordham \cite{blake} 2.5.1]
\label{thm:blake} Given a word $w \in F$ described by the reduced tree
pair diagram $(T_-,T_+)$, the length $|w|_{\F}$ of the word with respect
to the generating set $\F$ is the sum of
the weights of the caret pairings in $(T_-,T_+)$.
\end{theorem}

\subsection{How generators affect a tree pair diagram}

\label{sec:action}

The strength of Fordham's method is that it requires only
geometric information about the pair of trees representing an element
$w$ to determine $|w|$.  Beginning with an element $w=(T_-,T_+)$,
if we knew the reduced pair of trees which represented $w \alpha$ for
$\alpha \in \{ \x^{\pm 1}, \ x_1^{\pm 1} \}$, we could deduce the word length
of $w \alpha$.  We now discuss how the tree pair diagrams for $w$ and $w \alpha$ are related.

We begin with a lemma from Fordham \cite{blake} which states under fairly broad conditions,
that when applying a generator to a tree pair $(T_-,T_+)$
exactly one pair of caret types will change.

\begin{lemma}[Fordham \cite{blake}, Lemma 2.3.1]
\label{lemma:conditions}
Let $(T_-,T_+)$ be a reduced pair of trees, each having $m+1$
carets, representing an element $w \in F$, and $\alpha$ any
generator of $\F$.
\begin{enumerate}
\item
If $\alpha = \x$, we require that the left subtree of the root of $T_-$
is nonempty.

\item
If $\alpha = \inv$, we require that the right subtree of the root of $T_-$
is nonempty.

\item
If $\alpha = x_1$, we require that the left subtree of the right child
of the root of $T_-$ is nonempty.

\item
If $\alpha = x_1^{-1}$, we require that the right subtree of the right child
of the root of $T_-$ is nonempty.

\end{enumerate}
If the reduced tree pair diagram for $x \alpha$ also
has $m+1$ carets,
then there is exactly one $i$ with $0 \leq i \leq
m$ so that the pair of caret types of caret $i$ changes when
$\alpha$ is applied to $x$.
\end{lemma}

We now begin to understand geometrically the action of a generator
of $\F$ on a reduced tree pair $(T_-,T_+)$, and the corresponding change in
normal form. In this section we will assume that the conditions of lemma
\ref{lemma:conditions} are met by the generic elements with which we begin.
The following geometric lemma describing the action of the generators in
$\F$ on an element $w$ is proven in \cite{ct}.
Let $C_R$ denote the caret which is the right child of the root caret of
$T_-$, and $C_L$ the left child of the root.

\begin{figure}\includegraphics[width=3in]{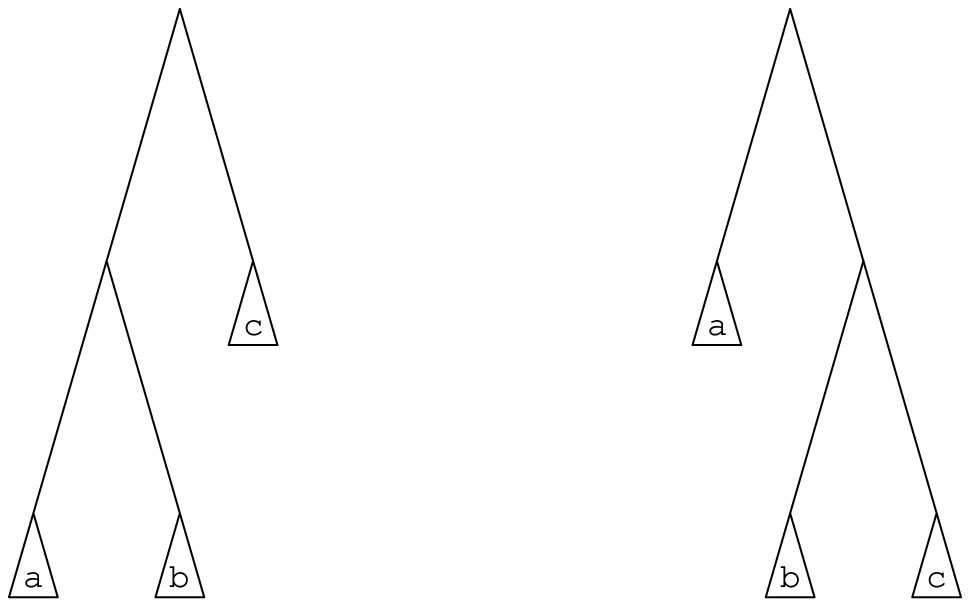}\\
\caption{The action of $\x$ transforms the right tree to the left
one, while the action of $\inv$ transforms the left tree to the
right one, in the figure above. Each tree represents only the negative
trees in their respective tree pairs. \label{x0invaction}}
\end{figure}

\begin{figure}\includegraphics[width=3in]{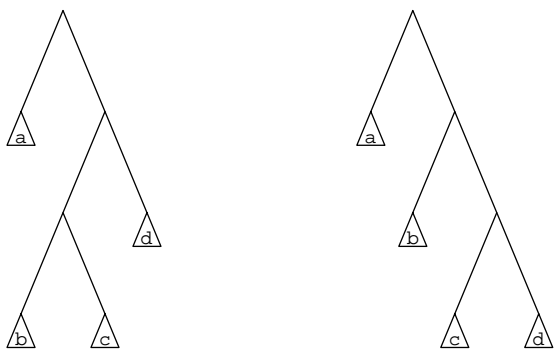}\\
\caption{The action of $x_1$ transforms the left tree to the right
one, while the action of $x_1^{-1}$ transforms the right tree to
the left one, in the figure above. Again, each tree represents only the
negative trees in their respective tree pairs.
\label{x1invaction}}
\end{figure}

\begin{lemma}[\cite{ct}, Lemmas 2.6,2.7]
\label{lemma:geoaction}
If $w=(T_-,T_+) \in F$ satisfies the appropriate condition of lemma \ref{lemma:conditions},
then $\x$ (resp. $\inv$) alters the position of the right subtree of
$C_L$ in $T_-$ (resp. the left subtree of $C_R$) as depicted in figure
\ref{x0invaction}.  In addition, $x_1$ and $x_1^{-1}$ perform
analogous operations on the subtrees of $C_R$, as depicted in figure
\ref{x1invaction}.
\end{lemma}

Notice that in all of the descriptions above, the tree $T_+$ from
the pair $w=(T_-,T_+)$ is not affected when a generator is applied to $w$.
This is not true in general for reduced tree pair diagrams not satisfying the
conditions of lemma \ref{lemma:conditions}.
In general, $T_+$ can be
affected in exactly three ways:

\begin{enumerate}
\item when $T_-$ has a single left edge, and the generator is
$\x$,

\item
when the left  subtree of the right child of the root caret of $T_-$ is empty,
and the generator is $x_1$, or

\item
if the generator is $\alpha$ and the pair of trees corresponding
to $x \alpha$ is not reduced.
\end{enumerate}

When the generators $\x$ and $\inv$ are applied to an element $w
\in F$, the change in normal form is straightforward.  Namely, $w
\inv$ remains in normal form.  If $w=w' \inv$ in normal form, then
$w \x = w'$ in normal form. Otherwise, $w = x_0^m w''$ in normal
form, where $m \geq 0$. In this case, $w \x = \x^{m+1} \phi(w'')$,
where $\phi:F \rightarrow F$ is the shift map which increases the
index of each generator in the normal form of $w$.

We now determine the change in normal form when a generator $x_1^{\pm 1}$ is
applied to an element $w$ in normal form.  The following lemmas are
proven in \cite{ct}.

\begin{lemma}[The normal form of $wx_1^{-1}$, \cite{ct}, Lemma 2.4]
\label{lemma:leafnumber1} Let $w  \in F$ be represented by the
tree pair $(T_-,T_+)$, and have normal form $x_{i_1}^{r_1} \cdots
x_{i_n}^{r_n}x_{j_m}^{-s_m} \cdots x_{j_1}^{-s_1}$.  Then $w
x_1^{-1}$ has normal form
\begin{equation}
\label{eqn:nf1}
 x_{i_1}^{r_1} \cdots x_{i_n}^{r_n}x_{j_m}^{-s_m}
\cdots x_{j_{q+1}}^{-s_{q+1}} x_{\alpha}^{-1} x_{j_q}^{-s_q}
\cdots x_{j_1}^{-s_1},
\end{equation}
where we might have $\alpha  =  j_{q+1}$.  If the root caret of $T_-$
has right and left subtrees $S_R$ and $S_L$ respectively, then
$\alpha$ is smallest leaf number in $S_R$.
\end{lemma}

\begin{lemma}[The normal form of $wx_1$, \cite{ct}, Lemma 2.5]
\label{lemma:leafnumber2} Let $w$ satisfy the conditions of lemma
\ref{lemma:conditions} and have normal form $x_1^{r_1} \cdots
x_{i_n}^{r_n}x_{j_m}^{-s_m} \cdots x_{j_1}^{-s_1}$.  Then $w x_1$
has normal form:
\begin{equation}
\label{eqn:nf2} x_1^{r_1} \cdots x_{i_n}^{r_n}x_{j_m}^{-s_m}
\cdots x_{j_l}^{-(s_l-1)} \cdots x_{j_1}^{-s_1},
\end{equation}
for some index $j_l$, which is the smallest
leaf number in the right subtree of $T_-$.

\end{lemma}

\subsection{Calculating distance between elements using tree pair diagrams}
\label{sec:distance}

When viewing elements of $F$ as homeomorphisms of $[0,1]$ it is clear that
inversion and group multiplication correspond to inversion and
composition of homeomorphisms.
We now interpret inversion and group multiplication in terms of
tree pair diagrams.

Inversion of a group element $f$ given by a tree pair diagram $(T_-,T_+)$
is simply the tree pair diagram $(T_+,T_-)$.   This can be seen easily from the
normal form and is interpreted via homeomorphisms in \cite{cfp}.

\begin{figure}\includegraphics[width=4.5in]{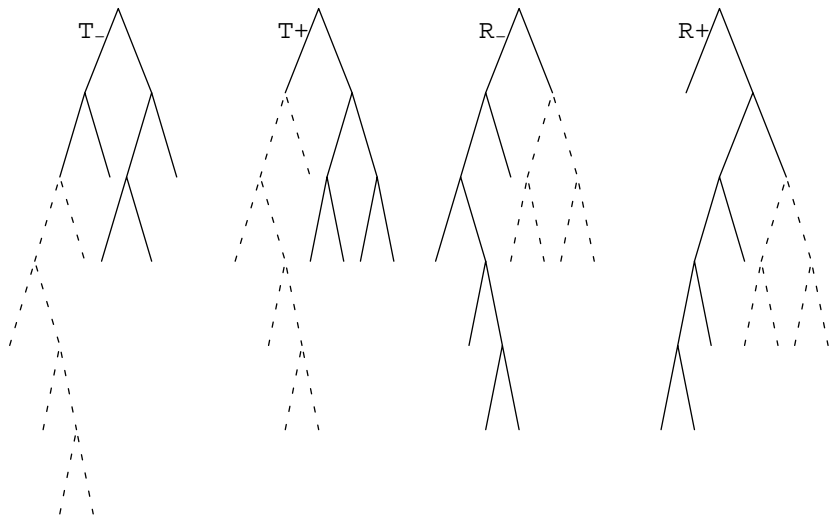}\\
\caption{Composing $x_1^3 x_2^{-1} x_1^{-1} x_0^{-2}$ with
$x_1 x_2^{-1} x_0^{-1}$ to
get $x_1^3 x_5 x_6^{-1} x_2^{-1} x_1^{-1} x_0^{-3} $ by adding the dashed carets.
\label{compexample}}
\end{figure}

Given two group elements $f,g \in F$ with tree pair diagrams
$f=(T_-,T_+)$ and $g=(R_-,R_+)$ we would like to form their product
$f g$ by a process consistent with composition of homeomorphisms.
When $T_+$ and $R_-$ are identical, we see that $g \circ f$ is represented
by the (possibly unreduced) tree pair diagram $(T_-,R_+)$.  This corresponds to composition of the
piecewise linear homeomorphisms represented by $f$ and $g$, where $Range(f)
=Domain(g)$.

When $T_+$ and $R_-$ differ, we create temporary, unreduced
representatives of $f$ and $g$ in which the new trees $T_+$ and $R_-$ are identical.
Then the composition is carried out in the same manner as described above.

Figure \ref{compexample}  gives an example of the composition of
two elements of $F$.  The solid lines indicate the original carets and the
dashed lines indicate carets added to perform the composition which create unreduced
representatives of both elements.  To measure the distance between two
elements of $F$, we use the word metric on the product $f^{-1}g$ to obtain the metric
$d(f,g)=|f^{-1}g|$.

\section{Estimating the word metric $d_{\F}$}
\label{sec:metric}

It follows immediately from the chart in \S
\ref{sec:Fordham} that the word length $|w|_{\F}$ of $w = (T_-,T_+) \in
F$ can be estimated in terms of the number of carets $N(w)$ in
either tree.  This estimate is analogous to the one given by Burillo,
Cleary and Stein
\cite{bcs} which has multiplicative constant $12$ for the upper bound.

\begin{theorem}
\label{thm:blake'sbound} Let $w \in F$ be represented by a tree
pair $(T_-,T_+)$ in which each tree has $N(w)$ carets. Then
$$N(w) - 2 \leq |w|_{\F} \leq 4N(w)-4.$$
\end{theorem}

\begin{proof}
We first note that every reduced tree pair has a caret type pair
$(L_0,L_0)$ of weight $0$.  Also, it is possible to have
the last caret type pair be $(R_0,R_0)$ which also has weight $0$.
(The only instance in which this does not happen is when the root
caret of $T_-$ or $T_+$ has no right subtree.)
For the upper bound, we can ignore the first $(L_0,L_0)$ caret pair,
and looking at the
chart in \S \ref{sec:Fordham}, we see that the maximum weight of
any other pair of caret types is $4$.
Thus the word length of $w$ is at most $4(N(w)-1)$.

To compute the lower bound, we ignore pairs of carets of type $(L_0,L_0)$
and $(R_0,R_0)$.  Any other pair of caret types has weight at least $1$,
and the lower bound is easily obtained.
\end{proof}

It is natural to ask if the multiplicative coefficient of 4  in theorem
\ref{thm:blake'sbound} can be improved to $3$, since
looking at the chart in \S \ref{sec:Fordham} we see that there are very few entries
which are $4$; that is, very few caret type pairs actually have weight $4$.  The
answer is no; one can produce words which get extremely close to
the bound of $4$ by pairing $I_R$ and $I_0$ carets in a particular
way.

\begin{figure}\includegraphics[width=4.5in]{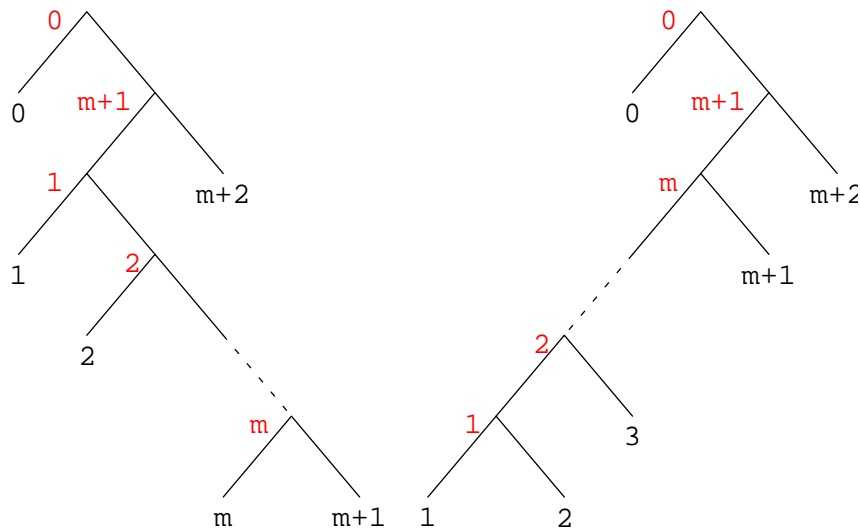}\\
\caption{ Tree pair diagram $(T_-,T_+)$ for the word
$x_1^m x_m^{-1} x_{m-1}^{-1} \cdots x_1^{-1}$ with the carets and
the leaves numbered.}
\label{iri0example}
\end{figure}

\begin{example} Words of the form
$x_1^m x_m^{-1} x_{m-1}^{-1} \cdots x_1^{-1}$ where $m>1$ is a
positive integer realize the upper bound of $4$ in theorem \ref{thm:blake'sbound}.
\end{example}

\noindent
Words of the above form are represented by the tree pair diagram given in
figure \ref{iri0example}.  They are constructed so that most carets
in $T_-$ are of type $I_R$ and are paired with carets of type $I_0$ in
$T_+$, to give a weight of $4$ per pair for most caret pairs. The
weights of the different caret pairs are summarized in the following table.

\smallskip
\begin{center}
\begin{tabular}{|c|c|c|c|}
\hline
Caret numbers & Caret types & Weight per caret & Total weight \\
\hline 0 & $(L_0,L_0)$ & 0 & 0 \\
$ 1, 2, \cdots ,m-1$ & $(I_R,I_0)$ & 4 & $4(m-1)$ \\
m & $(I_0,I_0)$ & 2 & 2 \\
m+1 & $(R_0,R_0)$ & 0 & 0  \\\hline
\end{tabular}
\end{center}
\smallskip

\noindent The total weight of a word $w$ of this form is $4(m-1) + 2 =
4m-2$. The total number of carets  $N(w)$ is $m+2$, so these
examples, which have weight $4m-2 = 4N(w)-10$, show that
for large $N(w)$ the multiplicative  coefficient of 4 in Theorem
\ref{thm:blake'sbound} is optimal.

It is also natural to wonder if the lower bound can be realized;
that is, are there examples of words $w$ where the number of carets is
exactly two more than the word length $|w |_{\F}$?  This will
always be true for words of the form  $x_1^{\pm
n}$ (but not for $\x^{\pm n}$).
In the following example, we show that it can also be true for
more complicated words.

\begin{figure}\includegraphics[width=4.5in]{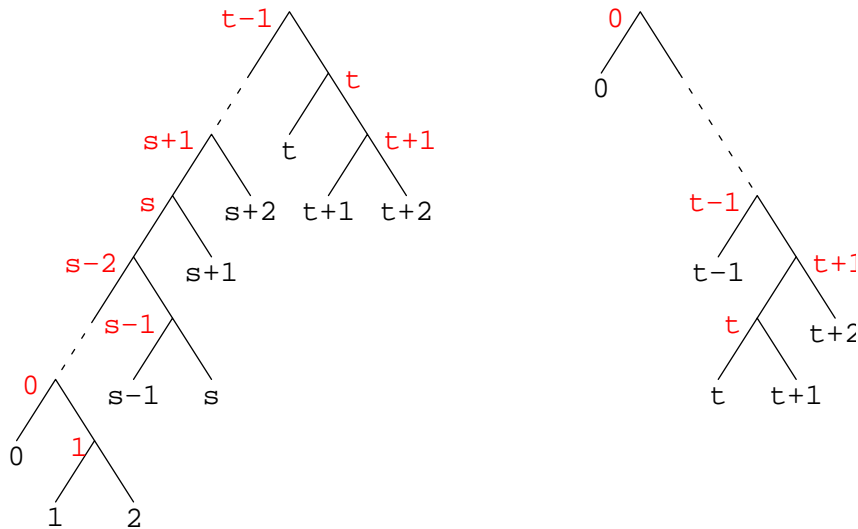}\\
\caption{The tree pair diagram $(T_-,T_+)$ for words
$x_t x_s^{-1} x_{s-2}^{-1} \cdots x_5^{-1} x_3^{-1} x_1^{-1}
x_0^{-m}$, with both the carets and the leaves numbered.}
\label{shortexample}
\end{figure}

\begin{example} Words of the form
$x_t x_s^{-1} x_{s-2}^{-1} \cdots x_5^{-1} x_3^{-1} x_1^{-1}
x_0^{-m}$ where $t> s+2$, $s$ is odd, and $m$ is chosen so that
the root caret of $T_-$ is caret number $t-1$ realize the lower bound
of theorem \ref{thm:blake'sbound}.  \end{example}
\noindent These words are
represented by tree pairs of the form given in Figure \ref{shortexample}.  The
weights of the caret type pairs are summarized in the following
table.

\smallskip
\begin{center}
\begin{tabular}{|c|c|c|c|}
\hline
Caret numbers & Caret types & Weight per caret\\
\hline 0 & $(L_0,L_0)$ & 0\\
$ 2, 4, \cdots ,s-1$, even numbers & $(L_l, \rni)$ & 1 \\
$1,3, \cdots ,s$, odd numbers & $(I_0,\rni)$ & 1\\
$s+1, \cdots t-2$ & $(L_L,\rni)$ & 1  \\
$t-1$ & $(\Ll, R_{I})$ & 1 \\
$t$ & $(R_0,I_0)$ & 1 \\
$t+1$ & $(R_0,R_0)$ & 0 \\
\hline
\end{tabular}
\end{center}
\smallskip

\noindent It is clear from the table that the total weight of the
word is two less than the number of carets, realizing the lower
bound of theorem \ref{thm:blake'sbound}.

Theorem \ref{thm:blake'sbound} has an immediate improvement for
strictly positive or negative words.

\begin{corollary}
Let $w$ be a strictly positive or negative word, represented by a
tree $T$ having $N(w)$ carets.  Then
$$N(w) - 2 \leq |w|_{\F} \leq 3N(w)-3.$$
\end{corollary}

\begin{proof}
Since $w$ is strictly positive or negative, one of the trees in
the tree pair diagram for $w$ consists entirely of $R_0$ carets
(excepting the root caret which is of type $L_0$ and does not contribute
to the total weight of the word.)
Thus, we only need to look at the first column of the chart in \S
\ref{sec:Fordham} to assign weights to the different carets.  We
notice that the maximum weight in the first column of this chart
is $3$, and the corollary follows.
\end{proof}

\section{Dead end elements}
\label{sec:deadend}

We now consider the Cayley graph $\Gamma$ of $F$ with respect to the finite generating
set $\F$.  Fordham describes a family of elements $w \in F$ which we call {\em dead end
  elements} that have the property that all four generators $x_0^{\pm
  1}$ and $x_1^{\pm 1}$ decrease the word length of $w$.  They are ``dead ends'' in
the sense that a geodesic ray in $\Gamma$ from the identity cannot be pass through them,
that is, a ray through a dead end element $w$ can no longer be geodesic past $w$.
We prove that all dead end elements have a particular form, which is
  slightly more general than the examples given by Fordham \cite{blake}, and
discuss other possible dead end behavior in $\Gamma$.

 \subsection{The form of dead end elements in $F$}

\begin{figure}\includegraphics[width=3in]{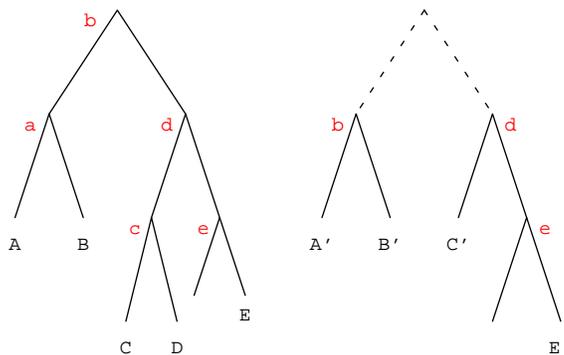}\\
\caption{The general form of dead-end elements of $F$.  Capital letters represent
(possibly empty) subtrees, while lower case letters label the carets.
Since E and E' are non-empty, caret d is of type $\rni$ in both trees.
\label{deadend}}
\end{figure}

\begin{theorem}
\label{thm:deadend}
All dead end elements in $F$ are given by tree pairs of the form
in Figure \ref{deadend} where the subtrees $E, A'$ and $E'$ are nonempty.
\end{theorem}

The key step in the proof of theorem \ref{thm:deadend} is enumerating the
conditions under which a specific generator $\alpha \in \{ \x^{\pm 1}, \ x_1^{\pm 1} \}$
decreases the word length of an element.
This requires a detailed understanding of caret pairings as well as
how a generator can affect a tree pair.  We begin with a lemma of Fordham
\cite{blake} which states that if $w$ does not satisfy the appropriate
condition of lemma \ref{lemma:conditions}, then
$|w \alpha| = |w| + 1$. This allows us to  consider only words satisfying
all the conditions of lemma \ref{lemma:conditions} as possible dead end words.

\begin{lemma}[Fordham \cite{blake}, Lemma 2.4.2]
\label{lemma:wrongform}
Let $\alpha$ be a generator of $\F$, and $w=(T_-,T_+)$ a word not satisfying
the condition in lemma \ref{lemma:conditions} corresponding to $\alpha$.
Then $|w \alpha| = |w| + 1$.
\end{lemma}

Thus, if $w \in F$ does not satisfy all the conditions of lemma \ref{lemma:conditions},
then $w$ is not a candidate for a dead end word.

We let $C_R$ and $C_L$ denote the left and right carets, respectively,
of the root caret of $T_-$.  Similarly, let $C_{RR}$ and $C_{RL}$
denote the right and left children, respectively, of the caret $C_R$.
Let $S_L$ and $S_R$ denote the left and right subtrees,
respectively, of the root caret of $T_-$.  Continue this notation to let
$S_{RL}$ denote the left subtree of $C_R$, etc.  In general, if
$\overline{D}$ is a string of entries from the set $\{ R, \ L \}$,
then $S_{\overline{D}L}$ is the left subtree of $C_{\overline{D}}$.
We analogously define $S_{\overline{D}R}$.
We use the notation $R_*$ to refer to a right caret of any type, and
$I_*$ to refer to an interior caret.

We now rewrite the chart in \S \ref{sec:Fordham} from a different perspective.
Assume that $w = (T_-,T_+)$ satisfies the conditions of lemma \ref{lemma:conditions}.
We are interested in the conditions on the pairings of certain carets in $T_-$
which determine whether $|w \alpha| = |x| + 1$ or $|w \alpha| = |x| - 1$
for a given generator $\alpha \in \{ \x^{\pm 1},  x_1^{\pm 1} \}$.

These conditions are summarized in the charts below.
Since we are only considering elements $w=(T_-,T_+) \in F$ satisfying
lemma \ref{lemma:conditions}, we know that only the type of a single caret $C$ in $T_-$
will change when a generator of $\F$ is applied.  We list the initial
type of this caret $C$ in the second column of the charts below, and the
new type of caret $C$ in the third column.  Column $4$, titled ``Increase", lists the types of
carets in $T_+$ which can be paired with $C$ in order for
$|w \alpha| = |w| + 1$, and in column $5$, titled ``Decrease", we put the pairings of $C$
which yield $|w \alpha| = |w| - 1$.
These pairings are determined by whether certain subtrees of $T_-$ are
empty or not.  These conditions are summarized in the following table.

\noindent
{\bf Changes in word length when a specific generator is applied to $w = (T_-,T_+)$.}

\smallskip
\noindent Consider the elements $w=(T_-,T_+)$ and $wx_0$. Caret
$C$ is the root caret of $T_-$.

\smallskip
\begin{tabular}{|c|c|c|c|c|}
\hline
  Condition & Initial type & New type & Increase if & Decrease if\\

on $T_-$ & of caret $C$ & of caret $C$& $C$ paired with & $C$ paired with \\
\hline
$S_{RL} \neq \emptyset$ & $\Ll$ & $ R_I $ & $R_*, I_* $ & $\Ll$ \\
\hline
$S_{RL} = \emptyset , S_{RR} \neq \emptyset$ & $ \Ll$ & $\rni $ & $ R_*, I_R $ & $\Ll, I_0$ \\
\hline
$S_{RL} = \emptyset, S_{RR} = \emptyset $ & $\Ll$  & $ R_0 $& $ \rni,R_I, I_R $&$ R_0, \Ll, I_0$
 \\
\hline

\end{tabular}

\bigskip
\noindent Consider the elements $w=(T_-,T_+)$ and $wx_0^{-1}$.
Caret $C$ is the caret $C_R$ of $T_-$.

\smallskip
\begin{tabular}{|c|c|c|c|c|}
\hline
  Condition & Initial type & New type & Increase if & Decrease if \\
 on $T_-$ & of caret $C$ & of caret $C$& $C$ paired with & $C$ paired with \\
\hline
$S_{RL} \neq \emptyset$ & $R_l$ &  $\Ll $& $\Ll$& $R_*, I_* $ \\
\hline
$S_{RL} = \emptyset, S_{RR} \neq \emptyset $ & $\rni$ & $\Ll $  & $\Ll, I_0$ & $ R_*, I_R $\\
\hline
$S_{RL} = \emptyset, S_{RR} = \emptyset $ & $R_0$  & $ \Ll $&$ R_0, \Ll, I_0$& $ \rni,R_I, I_R $ \\
\hline

\end{tabular}

\bigskip
\noindent Consider the elements $w=(T_-,T_+)$ and $wx_1$. Caret
$C$ is the caret $C_{RL}$ of $T_-$.

\smallskip
\begin{tabular}{|c|c|c|c|c|}
\hline
  Condition & Initial type & New type & Increase if & Decrease if \\
 on $T_-$ & of caret $C$ & of caret $C$& $C$ paired with & $C$ paired with \\
\hline
$S_{RLR} \neq \emptyset$ & $I_R$ & $ R_I $&none & any \\
\hline
$S_{RLR} = \emptyset, S_{RR} \neq \emptyset $ & $ I_R$  & $ \rni $  & $R_0,\rni$ & $ \Ll,I_*,R_I $\\
\hline
$S_{RLR} = \emptyset, S_{RR} = \emptyset $ & $ I_0$  & $ R_0$  & $\rni$ & $ \Ll,I_*,R_I ,R_0$\\
\hline
\end{tabular}

\bigskip
\noindent Consider the elements $w=(T_-,T_+)$ and $wx_1^{-1}$.
Caret $C$ is the caret $C_{R}$ of $T_-$.

\smallskip
\begin{tabular}{|c|c|c|c|c|}
\hline
  Condition & Initial type & New type & Increase if & Decrease if \\
 on $T_-$ & of caret $C$ & of caret $C$& $C$ paired with & $C$ paired with \\
\hline
$S_{RRL} \neq \emptyset$ & $R_I$ & $ I_R $& any & none  \\
\hline
$S_{RRL} = \emptyset, S_{RRR} \neq \emptyset $ & $\rni$  & $ I_0 $  & $ \Ll,I_*,R_I $ & $ \Ll,I_*,R_I $ \\
\hline
$S_{RRL} = \emptyset, S_{RRR} = \emptyset $ & $R_0$ & $ I_0$ & $ \Ll,I_*,R_I ,R_0$  & $\rni$  \\
\hline
\end{tabular}

\bigskip

\bigskip
\noindent
{\it Proof of theorem \ref{thm:deadend}.}
From lemma \ref{lemma:wrongform} we may assume
that the initial word $w$ satisfies
all the conditions of lemma \ref{lemma:conditions}.
Combining the above four charts, we see that for all four generators
to reduce the word length of $w$, the element $w$ must be represented
by a tree pair $(T_-,T_+)$, where $T_-$ is given in figure \ref{deadend}.
We now determine which of the following subtrees $A, \ B, \ C, \ D, \ $ and $E$
may be empty.

The possible pairings of carets $a,b,c$ and $d$ are also determined by
the above four charts.  The combination of these conditions restricts
the pairings further, as follows.
It is now clear that $\x$ causes $b$ to become a caret of type $R_I$, forcing
$b$ to be paired with a caret of type $L_L$.  Since $a<b$, caret $a$
must be paired with a caret appearing before $b$, thus the left
subtree $A'$ of caret $b$ in $T_+$ in nonempty.

When $x_1^{-1}$ reduces the word length of $w$, we now see that caret
$d$ becomes a caret of type $I_0$ rather than of type $I_R$, and thus
$d$ must be paired with a caret of type $R_0$ or $R_{NI}$.  It follows
that the left subtree of caret $d+1 = e$ in $T_+$ is empty (otherwise
caret $d$ in $T_+$ would be of type $R_I$).

We are now left with showing that the subtrees $E$ of $T_-$ and $E'$
of $T_+$ are nonempty.  Recall that the two trees $T_-$ and $T_+$ have
the same number of carets.  If $E$ was empty, then $E'$ must also be
empty, given the placement of caret $d$ in both trees.  If this is the
case, then the pair $(T_-,T_+)$ is not reduced, contradicting initial
assumptions.  Similarly, if $E'$ is empty, so is $E$ and the pair is
again not reduced.  Thus $E$ and $E'$ are both nonempty, and we have
shown that all dead end elements have the claimed form.
\qed

\subsection{Pockets in the Cayley graph of $F$}

A natural question to ask is whether there are more severe
forms of dead end phenomena in $F$.
Many researchers have wondered if there
might be {\em pockets} in the Cayley graph of $F$.
For $k>0$, the element $w \in F$ defines a $k$-pocket if $w \in B_{id}(n)$ where $n = |w|$, and
$B_w(k) \subset B_{id}(n)$; that is, if  all paths of length $k$ emanating from $w$
remain in the ball of radius $n$ centered at the identity.  A dead end element of
the form described above defines a $2$-pocket.  Though there are dead end elements,
we now show that there are not significantly deeper pockets.

\begin{theorem}There
are no $k$-pockets in the Cayley graph of $F$ with respect to the finite generating set $\F$
for $k \geq 3$.
\end{theorem}

\begin{proof}
A word $w$ which defines a $k$-pocket  must be a dead end word, otherwise
there would immediately be a path of length 1 from $w$ which left $B_{id}(n)$.
We will produce a path of length $3$ emanating from any dead end word $w$ which leaves $B_{id}(n)$.
The key fact in constructing this path is that according to theorem
\ref{thm:deadend}, in a dead end word $w$,  the left subtree of $C_{RR}$ is empty.
In figure \ref{deadend} we label the caret $C_{RR}$ is labelled $e$.
We label the exposed left leaf of caret $C_{RR}$ by $m$, and it follows that the caret
number of $C_{RR}$ is also $m$.

Let $w = (T_-,T_+) \in B_{id}(n)$ be a dead end word.  Then $|w \inv| = |w| - 1$ by construction.
In the tree pair diagram $(R_-,R_+)$ representing $w \inv$,  caret $m$ is the right child of
the root, and leaf $m$ is still its exposed left leaf.  We notice that $w \inv$ does not
satisfy the condition of lemma \ref{lemma:conditions} corresponding to $x_1$, so
an application of $x_1$ would increase the length of $w x_0$ by one to give $|w \inv x_1| = |w|=n$.
To construct this path of length 3 which leaves the ball, we look at  the resulting pair of trees for $w \inv x_1$.

If $w=(T_-,T_+)$ is a word not satisfying the condition of lemma
\ref{lemma:conditions} corresponding to the generator $\alpha$,
then $\alpha$ acts on $(T_-,T_+)$ in a way that adds an additional
caret.  To see this, we create an unreduced representative of $w$
by adding carets to both trees so that the unreduced
representative satisfies the condition of lemma
\ref{lemma:conditions} corresponding to $\alpha$.  Then, we allow
$\alpha$ to act on the unreduced pair; the resulting tree pair
will represent $w\alpha$ and be reduced.

An alternate way to determine the resulting pair of trees is to simplify
the normal form for $w x_{\alpha}$ and draw the corresponding trees.

Consider the tree pair $(R_-,R_+)$ representing the element $w \inv $.
By adding an extra caret attached to leaf $m$ to both trees in the pair,
we obtain an unreduced representative of $w \inv$ which satisfies the condition
of lemma \ref{lemma:conditions} corresponding to $x_1$.
Figure \ref{x1invaction} in \S \ref{sec:F} exhibits the change in
$R_-$ when $x_1$ is applied to this unreduced tree pair.  We see that
in the negative tree of the tree pair representing $w \inv x_1$, the right caret of the root
again has an exposed left leaf.
Hence the tree pair diagram of $w \inv x_1$ again does not satisfy the conditions
of lemma \ref{lemma:conditions} corresponding to $x_1$.
Thus, by lemma \ref{lemma:wrongform},  $|w \inv x_1^2| = |w| + 1 $ and is not in $B_{id}(n)$,
and there can be no $k$-pockets in the Cayley graph of $F$ for $k \geq 3$.
\end{proof}

\section{Counting carets}

Given the beautiful relationship between the normal form of elements
of $F$ and their representation as pairs of finite binary rooted trees,
it is natural to ask what information about the trees can be
readily determined from the normal form.  We show that the total number of
carets in each tree can be determined, as well as the number of
right, interior and left carets in each tree.
As an application, this count is used to give a more accurate
estimation of $|w|_{\F}$ than the one given in theorem \ref{thm:blake'sbound}.
We note that Burillo \cite{burillo} also estimates $|w|_{\F}$ from the normal
form of $w$.

We temporarily alter our notion of the normal form to make our computations
and formulae easier to
understand.  Namely, let $w$ have normal form
$$x_0^{r_0} x_{i_1}^{r_1} x_{i_2}^{r_2}\ldots x_{i_k}^{r_k} x_{j_l}^{-s_l}
\ldots x_{j_2}^{-s_2} x_{j_1}^{-s_1} x_0^{-s_0}$$ allowing for the possibility
that $s_0 = 0$ and $r_0 = 0$ if the generator $x_0$ does not appear in the
normal form.  We still retain the conditions for uniqueness of the normal form.

We first show that it is easy to detect from the normal form whether the
right subtree of the root caret of $T_{\pm}$ is empty.

\begin{lemma}[Seeing the right side of a tree]
\label{lemma:rightside}
Let the element $w = (T_-,T_+)$ have normal form
$x_0^{r_0} x_{i_1}^{r_1} x_{i_2}^{r_2}\ldots x_{i_k}^{r_k} x_{j_l}^{-s_l}
\ldots x_{j_2}^{-s_2} x_{j_1}^{-s_1} x_0^{-s_0}$.

\begin{enumerate}

\item
If $\displaystyle i_k < \Sigma_{m=0}^{k-1} r_m$  then the right subtree of
the root caret of $T_+$ is empty.

\item
If $\displaystyle j_l < \Sigma_{m=0}^{l-1} s_m$  then the right subtree of
the root caret of $T_-$ is empty.
\end{enumerate}
\end{lemma}

\begin{proof}
We work through the proof in the case of $T_-$. The proof is identical for
$T_+$.  We begin to build the tree $T_-$ using the fact that the exponent
of $x_j$ in the normal form is the leaf exponent of the leaf labelled $j$
in the tree.  Let
$S_L$ and $S_R$ denote, respectively, the left and right subtrees of the
root caret of $T_-$.

In the tree $T_-$, we build the subtrees $S_L$ and $S_R$, beginning with
$s_0 + 1$ left carets in $S_L$, with highest leaf number $s_0$.
If the index $j_1$ is greater than $s_0$, then we must add
$s_1$ interior carets added to the right subtree of the root of $T_-$.
If not, we add $s_1$ interior carets to the left subtree of the root of
$T_-$. Assume that these interior carets are added to the left subtree of
the root.  We then ask the question again, with the next index.
The highest leaf number in $S_L$ is now $s_0 + s_1$.
If $j_2 > s_0 + s_1$, then we must add
$j_2$ interior carets to the right subtree $S_R$, otherwise we add them to $S_L$.  As soon
as $j_n > \Sigma_{m = 1}^{n-1} s_m$ for some value of $n$, $S_R$ is
nonempty, and the analogous equation is still true for higher indices.
Thus it is sufficient to test the highest index to see if $S_R$ is empty.
\end{proof}

In the next proposition, we drop the requirement that the trees $T_-$
and $T_+$ in a tree pair have the same number of carets. In \S \ref{sec:F}
we saw that if two trees did not have the same number of carets,
this was easily achieved by adding extra
carets of type $R_0$ to one of the trees in the pair.  The addition of carets
had no affect on the normal form.  The next proposition
shows that one can compute the number of carets in each tree {\em which are not
these extra $R_0$ carets} directly from the normal form.

\begin{proposition}[Counting caret types from the normal form]
\label{prop:numberofcarets}
Let $w =(T_-,T_+)$ have normal form $x_0^{r_0} x_{i_1}^{r_1} x_{i_2}^{r_2}
\ldots x_{i_k}^{r_k} x_{j_l}^{-s_l}
\ldots x_{j_2}^{-s_2} x_{j_1}^{-s_1} x_0^{-s_0}$, where we do not
require that $T_-$ and $T_+$ have the same number of carets.

\begin{enumerate}
\item
If the subtree $S_R$ of $T_-$ is not empty, i.e. $w$ does not satisfy condition $(2)$ of
lemma \ref{lemma:rightside},
then the tree $T_-$ then has

\begin{enumerate}
\item $j_l + s_l + 1$ total carets,
\item
$s_0 + 1$ left carets,
\item
$\Sigma_{m=1}^l s_m$ interior carets, and
\item
$j_l + s_l -  \Sigma_{m=0}^l s_m $ right carets.
\end{enumerate}

\item
If the subtree $S_R$ of $T_-$ is empty, i.e. $w$ satisfies condition $(2)$ of lemma
\ref{lemma:rightside},
then the tree $T_-$ has

\begin{enumerate}
\item $\Sigma_{m=0}^l s_m + 1$ total carets,
\item
$s_0 + 1$ left carets,
\item
$\Sigma_{m=1}^l s_m$ interior carets, and
\item
no right carets.
\end{enumerate}

\end{enumerate}
\end{proposition}

\begin{proof}
We first prove case $(1)$, using the fact that the right
subtree of $T_-$ is not empty.
To see that the total number of carets of $T_-$ is $j_l + s_l + 1$, notice
that there must be a left leaf in $T_-$ labelled $j_l$, but no higher numbered
left leaves. Since all the remaining leaves are right, we need
$s_l + 1$ of them so
that all carets are complete (that is, have two leaves).

It follows from the definition of leaf exponent in \S \ref{sec:F}
that there are $s_0 + 1$ left carets in $T_-$.  Every interior
caret contributes $1$ to the exponent of an exposed leaf numbered
$i$, for $i \neq 0$. Thus the number of interior carets is given
by $\Sigma_{m=1}^l s_m$. It then follows that the number of right
carets of $T_-$ is $(j_l + s_m) - \Sigma_{m=0}^l s_m$; that is,
the total number of carets less the number of left and interior
carets.

To prove case $(2)$, we follow the proof for case $(1)$, omitting
the right carets.
\end{proof}

Note that we have the identical theorem for the positive tree
$T_+$, replacing any instance of $j_l$ with $i_l$ and $s_m$ with
$r_m$, and using condition $(1)$ of lemma \ref{lemma:rightside}.

When we reinstate the requirement that the trees in the pair have the
same number of carets $N$, the number $N$ will be the maximum of the
number of carets obtained below for $T_-$ and $T_+$.

We can now improve the bound in theorem
\ref{thm:blake'sbound} slightly using proposition
\ref{prop:numberofcarets}.
Theorem \ref{thm:blake'sbound} is imprecise when it must
assume, for example, in the upper bound, that each caret has
weight $4$.  In reality there are only a few caret type pairs
which carry weight $4$.

The theorem below has the assumption that $T_-$ has at most as
many carets of type $R_0$ as $T_+$.  This is equivalent to saying
that if we do not require $T_-$ and $T_+$ to have the same number
of carets, then $T_-$ has more carets than $T_+$.  This assumption
allows us to use the estimates in proposition
\ref{prop:numberofcarets} for the number of carets of either tree
in the pair $(T_-,T_+)$ when we reinstate the requirement that the
two trees have the same number of carets.

\begin{theorem}
\label{thm:betterbound} Let $w= x_0^{r_0} x_{i_1}^{r_1} x_{i_2}^{r_2}\ldots
x_{i_k}^{r_k} x_{j_l}^{-s_l} \ldots x_{j_2}^{-s_2} x_{j_1}^{-s_1} x_0^{-s_0}$
be in normal form, and let $(T_-,T_+)$ be the tree pair diagram
for $w$. Without loss of generality we assume that $T_-$ has at most
as many carets of type $R_0$ as $T_+$.

\begin{enumerate}
\item
If the subtree $S_R$ of $T_-$ is nonempty, then
$$
 j_l + s_l-1 \leq |w|_{\F} \leq 3(j_l + s_l)+  \Sigma_{m=1}^l s_m +2s_0.
$$

\item
If the subtree $S_R$ of $T_-$ is empty then
$$
\Sigma_{m=0}^l s_m - 1 \leq |w|_{\F} \leq 2s_0 + 4\Sigma_{m=1}^l s_m.
$$
\end{enumerate}
\end{theorem}

\begin{proof}
We begin with the proof of case $(1)$. The lower bound is
immediate, since from the chart in \S \ref{sec:F}, we know that
each caret except the first and possibly the last, has weight at
least $1$.  To obtain the upper bound, we use the expressions in
proposition \ref{prop:numberofcarets}, part $(1)$, for the number
of each type of caret. Looking at the chart in \S
\ref{sec:Fordham}, we note the maximum weight for a caret of each
of the three types: right, left and interior. Namely, the maximum
weight for a right caret is $3$, for a left caret is $2$ and for
an interior caret is $4$.  The upper bound is then  $3( \# \text{
right carets} ) + 2 ( \# \text{ left carets}-1 ) + 4( \# \text{
interior carets}).$  We omit a single left caret since there will
always be a caret of type $L_0$ in a pair of the form $L_0,L_0)$,
which has weight $0$.

In case $(2)$, the proof is identical, substituting the caret counts from
part $(2)$ of proposition \ref{prop:numberofcarets}.
\end{proof}

We note that for strictly positive or negative words, this method
may improve the bound significantly.

\begin{corollary}
\label{cor:pos/neg}
Let $w = x_0^{s_0} x_{j_1}^{s_1} x_{j_2}^{s_2}\ldots x_{j_k}^{s_k}$
or $w=x_{j_k}^{-s_k} \ldots x_{j_2}^{-s_2} x_{j_1}^{-s_1} x_0^{-s_0}$
be a strictly positive or negative word.
\begin{enumerate}
\item
If $w$ is a negative word which does not satisfy condition $(2)$ of lemma
\ref{lemma:rightside}, or a positive word which does not satisfy condition
$(1)$ of lemma \ref{lemma:rightside}, then
$$ 2(j_k + s_k) - \Sigma_{m=0}^k s_m -2 \leq |w|_{\F} \leq 2(j_k + s_k) +
\Sigma_{m=0}^k s_m-2s_0.$$

\item
If $w$ is a negative word which satisfies condition $(2)$ of lemma
\ref{lemma:rightside}, or a positive word which satisfies condition
$(1)$ of lemma \ref{lemma:rightside}, then
$$ \Sigma_{m=0}^k s_m  \leq |w|_{\F} \leq s_0 + 3 \Sigma_{m=1}^k s_m .$$
\end{enumerate}
\end{corollary}

\begin{proof}
We work through the proof in the case that $w$ is a strictly negative word, so
$w=x_{j_k}^{-s_k} \ldots x_{j_2}^{-s_2} x_{j_1}^{-s_1} $. The case when $w$ is a
positive word is completely analogous.

Let $w = (T_-,T_+)$, where $T_+$ is a tree consisting entirely of
the root and $R_0$ carets.  The proof of the upper bound in either case is
identical to the proof in theorem \ref{thm:betterbound} in the analogous
case, except
that when we look up on the chart to see the maximum weight of
each type of caret we get different results.  Namely, the maximum
weight for a right caret, when paired with an $R_0$ caret, is $2$
rather than the $3$ obtained in the proof of theorem
\ref{thm:betterbound}. Similarly, the maximum weight of a left
caret when paired with an $R_0$ caret is $1$, and of an interior
caret is $3$.  Using these new values, we compute the upper bound
in the corollary.

To obtain the lower bound, in either case, we know that each caret
has weight at least $1$, but any right carets paired with a caret
of type $R_0$ must have weight $2$. The exception is a pair of
types $(R_0,R_0)$, which has weight $0$.  But there can be at most
one pair of these types in a reduced word, and we account for it
by subtracting $1$.  Thus $(\# \text{carets}) + (\# \text{ right
carets}) -2 = (i_k + r_k + 1) + (i_k + r_k + 1) - \Sigma_{m=1}^k
r_m-2 = 2(i_k + r_k) - \Sigma_{m=1}^k r_m$ is the lower bound in
case $(1)$, as desired.  In case $(2)$ we obtain $s_1 +
\Sigma_{m=1}^k r_m + 1 + 0 -2 = s_1 + \Sigma_{m=1}^k r_m -1$ as
desired.
\end{proof}

\section{Construction of short paths}

\label{sec:paths}

We now address the following question.  Given a word $w \in F$ written
in normal form, how do we find a minimal length representative of $w$
in the finite presentation $\F$, that is, a string representing $w$ which
contains only $\x^{\pm 1}$ and $x_1^{\pm 1}$, whose length is equal
to $|w|_{\F}$?

While Fordham's methods present a simple way to calculate the word length
of a word presented in the normal form,
he does not present an easy way of obtaining minimal length representatives
in the presentation ${\mathcal F}$.  Fordham's  algorithms, including those implemented
as LISP programs in  his thesis \cite{blake}, do produce a
minimal length representatives for a general word in normal form, but the
method requires substantial checking of cases, and is more suitable for
computer than human execution.

There is a simple method for constructing a (usually non-minimal)
representative for a word given in normal form in terms of
$x_0^{\pm 1}$ and $x_1^{\pm 1}$, which we call the {\em
replacement method}. We simply use the relators of ${\mathcal P}$
to replace each occurrence of $x_n^{\epsilon}$ in the normal form
of $w$ by $\x^{-(n-1)} x_1^{\epsilon} \x^{n-1}$.  It is easy to
see that this method gives a string of generators which is not
generally a minimal representative for the initial word.  In
general, this path differs from the minimal length representative
by no more than a factor of $4$.

Below, we present a simple method to
obtain a minimal length representative for a strictly positive or
negative word in $F$; that is, a word whose normal form consists
entirely of generators with positive or negative exponents.
We call this method the {\em nested traversal method}.
For a general word $w$, the nested traversal method could
be applied separately to the negative and positive parts of
$w$, but in general that would
produce a representative of $w$ which is
not minimal in the finite presentation.

\subsection{The nested traversal method}
Assume that $w = (T_-,*)$ is a tree pair diagram for a negative
word with $m$ carets. The tree $*$ consists almost entirely of
$m-1$ carets of type $R_0$; the only other caret in $*$ is the
root caret which is of type $L_0$.  Every leaf in the tree $*$ has
leaf exponent $0$. Viewing the identity as the non-reduced pair
$(*,*)$, we detail a sequence of generators which transform it to
the pair $(T_-,*)$.   We note that it is never necessary to apply
the generator $x_1$ as part of this process.  Also, the sequence
of generators produced by this method does not give (in general) a
unique representative of the pair $(T_-,*)$.

For each caret type $C$, the nested traversal method produces a
sequence of generators necessary to transform a cart of type $R_0$
into a caret of type $C$.  For each caret type, this can be
accomplished in a fixed number of steps, described below, which is
exactly the weight of the pair $(C, R_0)$ in the chart in \S
\ref{sec:Fordham}. Specifically, the following chart gives the
number of steps (applications of a generator) required to
transform an $R_0$ caret into a caret of type $C$.  Note that the
single $L_0$ caret will come from the single existing $L_0$ root
caret in $*$.

\smallskip
\begin{center}
\begin{tabular}{|c|c|c|c|c|c|c|c|}
\hline Type of caret $C$ & $L_0$ & $L_L$ & $I_0$ & $I_R$ & $R_0$ & $\rni$ &
$R_I$ \\ \hline Number of steps & 0 & 1 & 1 & 3 & 0 & 2 & 2 \\
\hline
\end{tabular}
\end{center}
\smallskip
Thus it is clear that the length of the path given by the nested traversal
method (described below), once it is proven that the method produces the
appropriate tree, will be the word length of $w = (T_-,*)$.

We now describe the nested traversal method, and in the theorem
below, prove that it yields the desired tree $T_-$.  The goal of
this method is to identify a sequence of generators of $\F$ which
transform a caret of type $R_0$ into a caret of type $C$. In \S
\ref{sec:Fordham} we saw that applying $\alpha \in \{ \x^{\pm 1},
\ x_1^{-1} \}$ to a word $w$ satisfying the appropriate condition
of lemma \ref{lemma:conditions} will either change the caret type
of the root caret or the right child of the root caret of $T_-$.

First, we note that to transform an $R_0$ caret which is the right
child of the root caret of $T_-$ into an $L_L$ caret, we only need
apply the generator $\inv$. Similarly, the generator $x_1^{-1}$
transforms an $R_0$ caret which is the right child of the root
caret of $T_-$ into an $I_0$ caret.  This behavior was exhibited
in figures \ref{x0invaction} and \ref{x1invaction}.

The simplest example of a caret of type $I_R$ is an interior caret
with a single right child of type $I_0$.  We do not view creating
the $I_0$ caret (or whatever the right subtree of the $I_R$ caret
may contain) as part of the transformation of the $R_0$ caret into
the $I_R$ caret.  Thus, to create an $I_R$ caret, we must do the
following:  we apply $\inv$, moving the $R_0$ caret, which begins
as the right child of the root caret of $T_-$ and which we denote
$R$, to the root position of $T_-$.  Then we create the $I_0$
caret which will be the right child of the completed $I_R$ caret
(or create the appropriate right subtree of the completed $I_R$
caret) using the string of generators specified by the nested
traversal method for those caret or carets.  We do not count these
steps toward the creation the $I_R$ caret, but rather as the steps
required to create the caret or carets in the right subtree of the
final $I_R$ caret. Finally, we apply the generator $\x$, which
moves $R$ back to the position of right child of the root caret,
and apply $x_1^{-1}$, which moves $R$ to an interior caret while
moving the newly created $I_0$ caret (or the carets in the right
subtree of the $I_R$ caret) to the right subtree of $R$.  Thus,
the sequence of generators needed to create the $I_R$ caret is
$\inv \cdots \x x_1^{-1}$ where the $\cdots$ represents a sequence
of generators which create the carets in the right subtree of the
final $I_R$ caret,  and must occur before the $I_R$ caret can be
completed.

The type of a right caret $C$ is determined by the left subtree of
the right caret ``beneath" $C$ on the right side of the tree, that
is, the next right caret further from root than $C$.  Keeping in
mind that only one of two positions in the tree can be affected by
$\alpha \in \{ \x^{\pm 1}, \ x_1^{-1} \} $, namely the root
position and the right child of the root, it is clear how to
transform an $R_0$ caret, denoted $R$, which begins as the right
child of the root caret of $T_-$, into a caret of type $R_{NI}$ or
$R_I$. Namely, we apply $\inv$ to move $R$ to the root position of
the tree. Then we use the nested traversal method to create the
subtree of the right caret beneath $R$, remembering that these
generators are all involved in the creation of {\em other} carets.
Finally, we apply $\x$, moving $R$ back to the right side of the
tree, now with its right subtree correct and thus the caret type
is now of the appropriate type.

It is the nesting of these sequences of generators which gives the method its
name.  To summarize this method, consider the following chart, which lists
the sequences of generators required to create each type of caret.

\smallskip
\begin{center}
\begin{tabular}{|c|c|}
\hline Caret type &  \parbox{1.8in}{\begin{center}Generators
involved in creation of caret \end{center}}\\ \hline $L_0$ & none
\\ \hline $L_L$ & $x_0^{-1}$ \\ \hline $I_0$ &$ x_1^{-1}$ \\
\hline $I_R$ & $x_0^{-1} \ldots x_0 x_1^{-1}$\\ \hline $R_0$ &
none \\ \hline $R_{NI}$ & $x_0^{-1} \ldots x_0$ \\ \hline $R_I$ &
$x_0^{-1} \ldots x_0$\\ \hline
\end{tabular}
\end{center}
\smallskip

To use this method to produce a minimal length representative of a
positive or negative word, we traverse the tree in infix order.
When we encounter the caret types of $L_0$, $L_L$, $I_0$ and
$R_0$, we record at most a single appropriate generator. When we
encounter the other caret types $I_R$, $R_{NI}$ and $R_I$, we
record the initial generator of the sequence given in the chart
above, then apply the process recursively on the right subtree (if
present) of the caret; after completion of the process on the
right subtree, we record the one or two additional generators
specified in the chart, and continue with the infix traversal of
the tree $T_-$.

We will prove the following theorem.

\begin{theorem}
\label{thm:nestedtraversal} Let $w$ be a strictly positive or negative
word.  Then the nested traversal method produces a minimal length
representative for $w$.
\end{theorem}

We begin with an example of the nested traversal method.

\begin{example}
\label{example:nestedtrav} To understand the traversal
construction method, we use it to construct a minimal length
representative of the strictly negative word $w$ with normal form
$x_{10}^{-1} x_7^{-1} x_6^{-1} x_4^{-1} x_2^{-2} x_0^{-2}$.
\end{example}

\noindent If $w = (T_-,*)$, where $w$ is the word given above,
then the carets of $T_-$, in infix order, have caret types given
in the following table.

\begin{center}
\begin{tabular}{|c|c|c|c|c|c|c|c|c|c|c|c|c|}
\hline Caret number & 0 & 1 & 2 & 3 & 4 & 5 & 6 & 7 & 8 &9 &10 &
11\\ \hline Caret type & $L_0$ & $L_L$& $I_0$ &$ I_R$& $I_0$ &
$L_L$ & $I_R$& $I_0$ & $R_{NI}$ & $R_I$& $I_0$ & $R_0$ \\ \hline
\end{tabular}
\end{center}
\smallskip

Since a tree pair diagram for $w$ will require 11 carets in each
tree, we begin with a tree $*$ consisting of one caret of type
$L_0$ and 10 carets of type $R_0$. Caret $0$ is of type $L_0$ and
is identical to the $L_0$ caret of $*$. To create caret $1$, we
apply $x_0^{-1}$. To create caret $2$, we apply $x_1^{-1}$. Caret
$3$ is more complicated.  We begin with $\inv$, and then must
traverse the right subtree of caret $3$, which in this case is the
single caret numbered $4$. Since caret $4$ is of type $I_0$, it
requires only the generator $x_1^{-1}$. We now return to complete
the sequence of generators necessary to create caret $3$, namely
$x_0 x_1^{-1}$.  Caret $5$, the root caret, is of type $L_L$ and
requires only the generator $\inv$.  Caret $6$ is again more
complicated, since it is of type $I_R$. We begin with an
$x_0^{-1}$ before descending into the right subtree of caret $6$,
i.e. caret $7$, which is created via an $x_1^{-1}$.  we then
finish the creation of caret $6$ with the string $x_0 x_1^{-1}$.
To create caret $8$, we apply $\inv$ and begin to traverse the
right subtree of caret $8$.  For caret $9$, we apply the generator
$\inv$ and begin to traverse the right subtree of caret $9$. Caret
$10$ is created with a single $x_1^{-1}$ generator, and caret $11$
is simply the last $R_0$ caret from the initial tree $*$ and
requires no additional generators.  We then return to finish the
creation of caret $9$ by applying $\x$ to make it once again a
right caret, and we apply another generator $x_0$ to finish
creating caret $8$.  Thus we have produced the tree $T_-$
representing $w$ and found the string
$$
\inv x_1^{-1} \inv  x_1^{-1} \x  x_1^{-1} \inv \inv  x_1^{-1} \x
x_1^{-1} \inv \inv  x_1^{-1} \x \x
$$
to be a minimal representative for $w$.

\begin{figure}\includegraphics[width=6in]{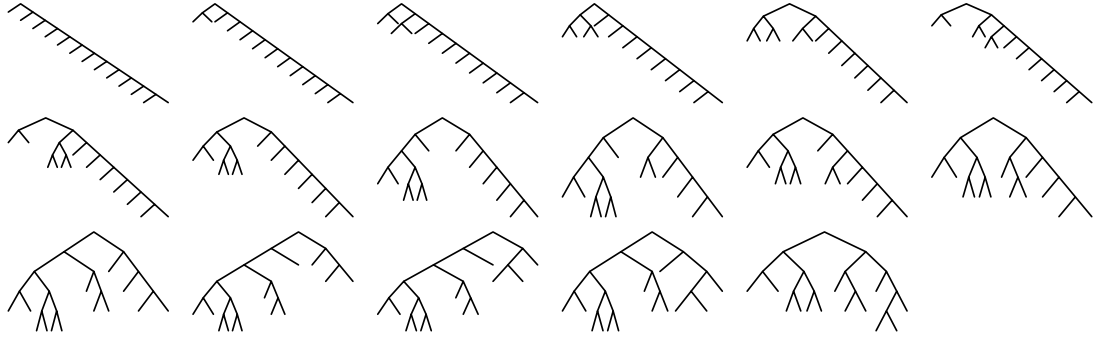}\\
\caption{The negative trees for the complete construction of a minimal length representative
for $w=x_{10}^{-1} x_7^{-1} x_6^{-1} x_4^{-1} x_2^{-2} x_0^{-2}$
via the nested traversal method.  The generators needed for each
step are exhibited below. The positive trees $*$ are omitted. \label{nestedtrav}}
\end{figure}

This entire process is illustrated in Figure \ref{nestedtrav}, and
described in more detail in the following four figures.

\begin{figure}\includegraphics[width=5in]{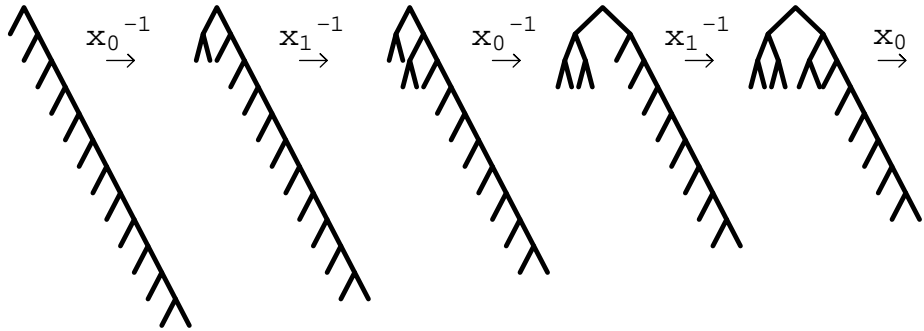}\\
\caption{The negative trees for the first four steps of  the nested traversal construction
of a minimal length representative for $w=x_{10}^{-1} x_7^{-1}
x_6^{-1} x_4^{-1} x_2^{-2} x_0^{-2}$, listing the generators
applied at each stage.  The positive trees $*$ are again omitted. \label{nestedtrav1}}
\end{figure}

\begin{figure}\includegraphics[width=5in]{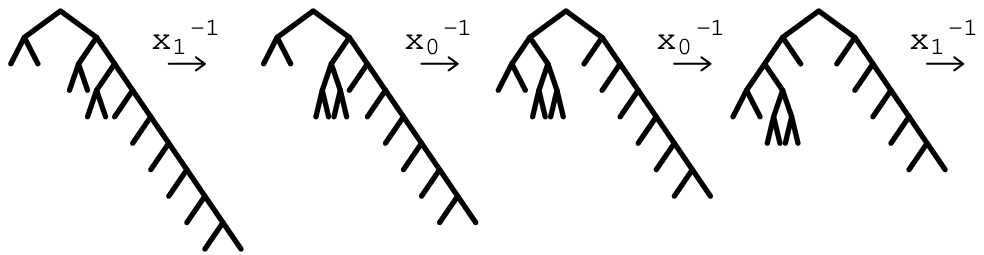}\\
\caption{The negative trees for the next four steps of the construction.
\label{nestedtrav2}}
\end{figure}

\begin{figure}\includegraphics[width=5in]{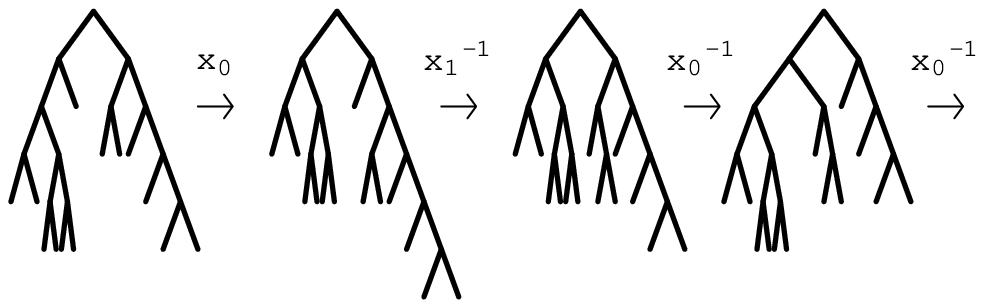}\\
\caption{The negative trees for the next four steps of the construction.
\label{nestedtrav3}}
\end{figure}

\begin{figure}\includegraphics[width=5in]{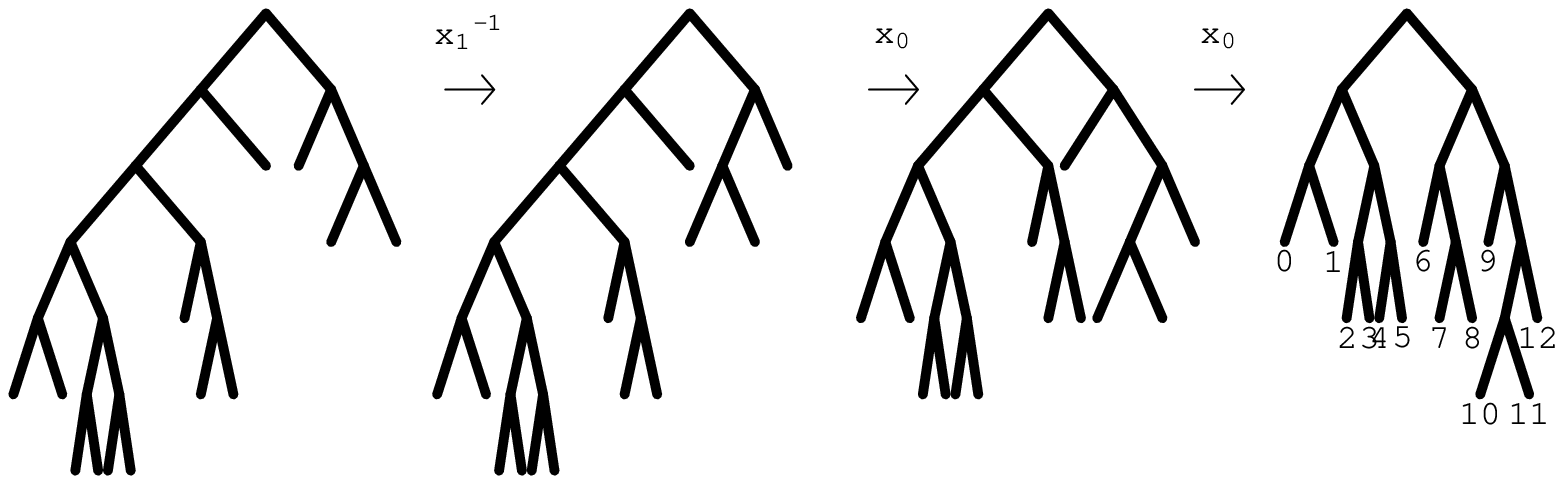}\\
\caption{The negative trees for the final four steps of the construction of a minimal
length representative for $w=x_{10}^{-1} x_7^{-1} x_6^{-1}
x_4^{-1} x_2^{-2} x_0^{-2}$ via the nested traversal method.  It
is easily checked that the final tree represents $w$.
\label{nestedtrav4}}
\end{figure}

We now prove theorem \ref{thm:nestedtraversal}.

\begin{proof}
We first prove that the nested traversal method produces the
correct tree for strictly negative words in $F$.  Once this is
proven, it is clear that the string of generators produced by the
method is a minimal length representative, because the number of
carets needed to create a caret of type $C$ is exactly the weight
of the pair $(C, R_0)$. We prove that the correct tree is produced
via induction on the length of the normal form of a negative word
$w \in F$.  Note that the length of the normal form of $w$ is
exactly $|w|_{\mathcal P}$, the word length in the infinite
presentation ${\mathcal P}$.

We begin with the base case of the induction, and $w = x_k^{-1}$.
 We must prove that the nested
traversal method produces the correct tree pair for $w$. If $w=
\inv$ then it is clear that the nested traversal method yields the
correct pair of trees.

If $w=x_k^{-1}$ with $k \neq 0$, it is easily checked that the tree $T_-$
for $w$ is of the following form: caret $0$ is of type $L_0$ and is
followed by $k-2$ carets of type $R_{NI}$, a single caret of type $R_I$,
a single caret of type $I_0$ and a final $R_0$ caret.  Creating these carets
in infix order, following the nested traversal method, we begin with
the tree $*$ containing a single $L_0$ caret, and $k+1$ carets of type $R_0$.

According to the nested traversal method, we apply $x_0^{-(k-1)}$,
which moves the first $k-1$ carets of type $R_0$ to the left side
of the tree.  Each of these generators is the first half of the
pair of generators needed for creating an $R_I$ or an $R_{NI}$
caret. Now, we create the $I_0$ caret by applying the generator
$x_1^{-1}$; this creates an $I_0$ caret with exposed leaves
numbered $k$ and $k+1$.  We finish with $x_0^k$, which moves all
but the initial $L_0$ caret on the left side of the tree back to
the right side of the tree. Thus, the nested traversal method,
with the correct nesting, produced the tree $T_-$ representing the
element $x_k^{-1}$.

We complete the proof by induction on the size of the negative
word.   We assume inductively
that the nested traversal method produces the correct tree for a
negative word $w=x_{j_l}^{-s_l} \ldots
x_{j_2}^{-s_2} x_{j_1}^{-s_1} $ with tree pair diagram $(T_-,*)$.
We consider a new word  $v=x_k^{-1} w $ in normal form whose length (in the infinite
presentation) is one more than the length of $w$. If
$v=(S_-,*)$, we prove by induction that $S_-$ can be constructed via the nested
traversal method.  We divide the proof into two cases:
\begin{enumerate}
\item
in case $1$, we assume that the normal form of $w$ begins with $x_k^{-1}$, and
\item
in case $2$, we assume  that the normal form of $w$ does not begin with $x_k^{-1}$.
\end{enumerate}
In both cases, the tree $S_-$ has at least one more caret than the tree $T_-$.

\medskip
\noindent
{\bf Case 1.} Let $w = x_k^{-m} w'$ in normal form, and thus $v=x_k^{-(m+1)} w'$ in
normal form.
In terms of tree pair diagrams, let $w = (T_-,*)$ and $v = (S_-,*)$.
Let $\gamma$ be the string of letters obtained via the nested
traversal method which creates the tree $T_-$, according to the induction
hypothesis.  We will obtain a string $\gamma'$ which is a minimal representative
of $v$ and creates the tree $S_-$.

\begin{enumerate}
\item
If $k=0$, it is clear that the nested traversal method produces the correct tree
$S_-$ and gives a minimal length representative for $v$.

\item Assume that $k \neq 0$.  If the right subtree of the root
caret of $T_-$ is empty (which we note can be detected via lemma
\ref{lemma:rightside}) then except for the initial $L_0$ caret,
$T_-$ consists entirely of left and interior carets. Since $k \neq
0$, there must be at least one interior caret in $T_-$.  Since
$\gamma$ constructs $T_-$ according to the nested traversal
method, we know that $\gamma$ can be written as $\gamma = \eta_1
x_1^{-1} \eta_2$, where the $x_1^{-1}$ creates the final interior
caret, which has caret number $k$.   In addition, $\eta_2$ is a
sequence of letters containing at least $l$ repetitions of $\inv$,
which create a series of left carets in $T_-$: the left carets
with caret numbers greater than $k$ as well as the left caret with
the highest caret number less than $k$.  The string $\eta_2$ may
contain other letters prescribed by the nested traversal method.

From the normal form of $w$, we know that the highest numbered
left leaf in $T_-$ is leaf number $k$, and it must be the exposed
left leaf of a string of $m$ carets of type $I_0$.  None of the
carets in this string can be of type $I_R$; if they were, then
there would be an exposed left leaf with positive leaf exponent
whose leaf number is greater than $k$.  Since $k$ is the largest
index of a generator in the normal form of $w$, such an $I_R$
caret cannot occur in this string of carets.

Such a string of $I_0$ carets is created according to the nested
traversal method via a string of $x_1^{-1}$ letters.  Thus we can
further enumerate $\gamma$ as $\gamma = \eta_1 x_1^{-m} \eta_2$.
Since the normal form of $v$ simply contains one more occurrence
of the generator $x_k$, we know that in $S_-$, the exposed leaf
numbered $k$ is the left exposed leaf of a sequence of $m+1$
carets of type $I_0$. Thus the tree $S_-$ is created by the string
of letters $\gamma' = \eta_1 x_1^{-(m+1)} \eta_2$. Since the
generator $x_1^{-1}$ is what the nested traversal method
determines will create an $I_0$ caret, we see that $\gamma'$ is
the string produced by the nested traversal method which generates
$S_-$.

\item Now suppose that the right caret of the root of $T_-$ is not
empty.  In this case, $\gamma$ can be written as $\gamma = \xi_1
x_1^{-m} \xi_2$, where $\xi_2$ accounts for the right carets that
are created in a nested manner.  It is easy to see that as above,
the tree $S_-$ is produced from the string $\gamma = \xi_1
x_1^{-(m+1)} \xi_2$, which is exactly the string of letters
produced by the nested traversal method to generate $S_-$.
\end{enumerate}

\medskip
\noindent {\bf Case 2.} We now assume that $v=x_k^{-1} w$ in
normal form, and that the generator $x_k$ does not appear in the
normal form of $w$.  In all of the subcases below, we begin with a
tree $*$ containing at least one more $R_0$ caret than the tree
used to create $T_-$.  Since the generator $x_k^{-1}$ does not
appear in the normal form of $w$, in the tree $T_-$, there is a
caret $C$ with an exposed leaf numbered $k$, which has leaf
exponent $0$.  It may be the case that one must add carets of type
$R_0$ to the right side of $T_-$ to obtain the caret $C$ with
exposed leaf $k$.  However, adding these carets does not affect
the normal form of the element. Again, let $\gamma$ be the string
of letters generated by the nested traversal method which creates
the tree $T_-$.

\begin{enumerate}
\item First, suppose that $C$ is an interior caret, which must be
of type $I_0$ and have an exposed right leaf numbered $k$ in order
for the leaf exponent of $k$ to be $0$. Since $k$ is larger than
any index of a generator in the normal form of $w$, it also
follows that the highest index appearing in the normal form of $w$
is at most $k-1$.

We can write $\gamma = \eta_1 x_1^{-1} \eta_2$, where $x_1^{-1}$ is the letter in
$\gamma$ which creates the caret $C$ from a right caret.  It is easily seen that
in $\gamma' = \eta_1 \inv x_1^{-1} \x x_1^{-1} \eta_2$, the caret $C$ becomes an
$I_R$ caret, the leaf exponent of $k$ is now $1$, and no new carets with higher
numbered leaves and positive leaf exponent are added, which might cause additional
generators to appear in the normal form for $v$.  We have exactly added the
sequence of generators corresponding to an $I_R$ caret in the chart above describing
the nested traversal method.

\item
Now suppose that $C$ is a left caret with an exposed right leaf numbered $k$.
Write $\gamma = \eta_1 \inv \eta_2$ where $\inv$ is the letter in $\gamma$ which
creates the caret $C$ according to the nested traversal method.  Then in the
tree created from $*$ by the string $\eta_1 x_1^{-1}$, the right child of the
root caret has a left subtree consisting of a single $I_0$ caret with exposed
leaves numbered $k$ and $k+1$.  Then, in the tree corresponding to the string
$\eta_1 x_1^{-1} \inv$, we have created the left caret $C$ whose right subtree
contains a single $I_0$ caret, with exposed leaves numbered $k$ and $k+1$.  Thus
the leaf exponent of $k$ is now $1$, and so $x_k^{-1}$ appears in the normal
form of the element.  Since the trees $T_-$ and $S_-$ differ only in this one place,
we see that the string $\gamma' = \eta_1 x_1^{-1} \inv \eta_2$ creates the tree
$S_-$.  The trees $S_-$ and $T_-$ differ in a single $I_0$ caret, and $\gamma'$
and $\gamma$ differ only in the generator $x_1^{-1}$, which is the letter that
the nested traversal method uses to create an $I_0$ caret, we see that $\gamma'$
is the string produced by this method for the tree $S_-$.

\item Finally, suppose that $C$ is a right caret in $T_-$ with
exposed left leaf $k$, with leaf exponent $0$.  The $C$ must be a
caret of type $R_0$; if it was not, then $k$ would be smaller than
the index of some generator appearing in the normal form of $w$,
contradicting initial assumptions.

Since $w$ is a negative word, it contains a single caret of type
$R_0$ which has a single right exposed leaf numbered $l$.  If
there was a caret of type $R_0$ with two exposed leaves,  the pair
$(T_-,*)$ could be reduced. It is necessary to add a string of
additional $R_0$ carets to the original tree $*$ in order to
obtain $S_-$.

Let $m=k-t$, where $t$ is the highest index of a generator appearing in the
normal form of $w$, so in particular $t<k$.  It now follows that the string
$\gamma'=\gamma x_0^{-m} x_1^{-1} x_0^m$ creates the tree $S_-$.  The first
$x_0^{-m}$ letters move the new $R_0$ carets to the left side of the tree,
the $x_1^{-1}$ creates an interior caret with exposed leaves labelled $k$ and
$k+1$, and the $x_0^m$ letters move the $m$ carets back to the right side
of the tree where they become of type $R_{NI}$ or $R_{I}$.  So again we see
that the additional letters needed exactly coincide with the letters
prescribed by the nested traversal method for constructing $S_-$.
\end{enumerate}

If we begin with a strictly positive word $w$, we use the nested
traversal method to construct a minimal length representative for
the strictly negative word $w^{-1}$.  The inverse of this path
will then produce a minimal length representative for $w$.
\end{proof}

The authors would like to thank Azer Akhmedov for helpful comments on an earlier
draft.

\bibliography{thompc}
\bibliographystyle{plain}

\bigskip

\begin{small}
\noindent Sean Cleary \\
Department of Mathematics \\
City College of New York \\
City University of New York \\
New York, NY 10031 \\
E-mail: {\tt cleary@sci.ccny.cuny.edu}

\medskip

\noindent
Jennifer Taback\\
Department of Mathematics and Statistics\\
University at Albany\\
Albany, NY 12222\\
E-mail: {\tt jtaback@math.albany.edu}
\end{small}

\end{document}